\magnification=1200

\hsize=11.25cm    
\vsize=18cm       
\parindent=12pt   \parskip=5pt     

\hoffset=.5cm   
\voffset=.8cm   

\pretolerance=500 \tolerance=1000  \brokenpenalty=5000

\catcode`\@=11

\font\eightrm=cmr8         \font\eighti=cmmi8
\font\eightsy=cmsy8        \font\eightbf=cmbx8
\font\eighttt=cmtt8        \font\eightit=cmti8
\font\eightsl=cmsl8        \font\sixrm=cmr6
\font\sixi=cmmi6           \font\sixsy=cmsy6
\font\sixbf=cmbx6

\font\tengoth=eufm10 
\font\eightgoth=eufm8  
\font\sevengoth=eufm7      
\font\sixgoth=eufm6        \font\fivegoth=eufm5

\skewchar\eighti='177 \skewchar\sixi='177
\skewchar\eightsy='60 \skewchar\sixsy='60

\newfam\gothfam           \newfam\bboardfam

\def\tenpoint{
  \textfont0=\tenrm \scriptfont0=\sevenrm \scriptscriptfont0=\fiverm
  \def\rm{\fam\z@\tenrm}
  \textfont1=\teni  \scriptfont1=\seveni  \scriptscriptfont1=\fivei
  \def\oldstyle{\fam\@ne\teni}\let\old=\oldstyle
  \textfont2=\tensy \scriptfont2=\sevensy \scriptscriptfont2=\fivesy
  \textfont\gothfam=\tengoth \scriptfont\gothfam=\sevengoth
  \scriptscriptfont\gothfam=\fivegoth
  \def\goth{\fam\gothfam\tengoth}
  
  \textfont\itfam=\tenit
  \def\it{\fam\itfam\tenit}
  \textfont\slfam=\tensl
  \def\sl{\fam\slfam\tensl}
  \textfont\bffam=\tenbf \scriptfont\bffam=\sevenbf
  \scriptscriptfont\bffam=\fivebf
  \def\bf{\fam\bffam\tenbf}
  \textfont\ttfam=\tentt
  \def\tt{\fam\ttfam\tentt}
  \abovedisplayskip=12pt plus 3pt minus 9pt
  \belowdisplayskip=\abovedisplayskip
  \abovedisplayshortskip=0pt plus 3pt
  \belowdisplayshortskip=4pt plus 3pt 
  \smallskipamount=3pt plus 1pt minus 1pt
  \medskipamount=6pt plus 2pt minus 2pt
  \bigskipamount=12pt plus 4pt minus 4pt
  \normalbaselineskip=12pt
  \setbox\strutbox=\hbox{\vrule height8.5pt depth3.5pt width0pt}
  \let\bigf@nt=\tenrm       \let\smallf@nt=\sevenrm
  \normalbaselines\rm}

\def\eightpoint{
  \textfont0=\eightrm \scriptfont0=\sixrm \scriptscriptfont0=\fiverm
  \def\rm{\fam\z@\eightrm}
  \textfont1=\eighti  \scriptfont1=\sixi  \scriptscriptfont1=\fivei
  \def\oldstyle{\fam\@ne\eighti}\let\old=\oldstyle
  \textfont2=\eightsy \scriptfont2=\sixsy \scriptscriptfont2=\fivesy
  \textfont\gothfam=\eightgoth \scriptfont\gothfam=\sixgoth
  \scriptscriptfont\gothfam=\fivegoth
  \def\goth{\fam\gothfam\eightgoth}
  
  \textfont\itfam=\eightit
  \def\it{\fam\itfam\eightit}
  \textfont\slfam=\eightsl
  \def\sl{\fam\slfam\eightsl}
  \textfont\bffam=\eightbf \scriptfont\bffam=\sixbf
  \scriptscriptfont\bffam=\fivebf
  \def\bf{\fam\bffam\eightbf}
  \textfont\ttfam=\eighttt
  \def\tt{\fam\ttfam\eighttt}
  \abovedisplayskip=9pt plus 3pt minus 9pt
  \belowdisplayskip=\abovedisplayskip
  \abovedisplayshortskip=0pt plus 3pt
  \belowdisplayshortskip=3pt plus 3pt 
  \smallskipamount=2pt plus 1pt minus 1pt
  \medskipamount=4pt plus 2pt minus 1pt
  \bigskipamount=9pt plus 3pt minus 3pt
  \normalbaselineskip=9pt
  \setbox\strutbox=\hbox{\vrule height7pt depth2pt width0pt}
  \let\bigf@nt=\eightrm     \let\smallf@nt=\sixrm
  \normalbaselines\rm}

\tenpoint

\def\pc#1{\bigf@nt#1\smallf@nt}         \def\pd#1 {{\pc#1} }

\catcode`\;=\active
\def;{\relax\ifhmode\ifdim\lastskip>\z@\unskip\fi
\kern\fontdimen2  -1.2 \fontdimen3 \string;}

\catcode`\:=\active
\def:{\relax\ifhmode\ifdim\lastskip>\z@\unskip\fi\penalty\@M\ \fi\string:}

\catcode`\!=\active
\def!{\relax\ifhmode\ifdim\lastskip>\z@
\unskip\fi\kern\fontdimen2  -1.1 \fontdimen3 \string!}

\catcode`\?=\active
\def?{\relax\ifhmode\ifdim\lastskip>\z@
\unskip\fi\kern\fontdimen2  -1.1 \fontdimen3 \string?}

% \catcode`\«=\active 
% \def«{\raise.4ex\hbox{%
%  $\scriptscriptstyle\langle\!\langle$}}
% 
% \catcode`\»=\active 
% \def»{\raise.4ex\hbox{%
%  $\scriptscriptstyle\rangle\!\rangle$}}

\frenchspacing

\def\raggedbottom{\topskip 10pt plus 36pt\r@ggedbottomtrue}

\def\pointir{\unskip . --- \ignorespaces}

\def\Medbreak{\vskip-\lastskip\medbreak}

\long\def\th#1 #2\enonce#3\endth{
   \Medbreak\noindent
   {\pc#1} {#2\unskip}\pointir{\it #3}\smallskip}

\def\decale#1{\smallbreak\hskip 28pt\llap{#1}\kern 5pt}
\def\decaledecale#1{\smallbreak\hskip 34pt\llap{#1}\kern 5pt}
\def\puce{\smallbreak\hskip 6pt{$\scriptstyle\bullet$}\kern 5pt}

\def\eqalign#1{\null\,\vcenter{\openup\jot\m@th\ialign{
\strut\hfil$\displaystyle{##}$&$\displaystyle{{}##}$\hfil
&&\quad\strut\hfil$\displaystyle{##}$&$\displaystyle{{}##}$\hfil
\crcr#1\crcr}}\,}

\catcode`\@=12

\showboxbreadth=-1  \showboxdepth=-1

\newcount\numerodesection \numerodesection=1
\def\section#1{\bigbreak
 {\bf\number\numerodesection.\ \ #1}\nobreak\medskip
 \advance\numerodesection by1}

\mathcode`A="7041 \mathcode`B="7042 \mathcode`C="7043 \mathcode`D="7044
\mathcode`E="7045 \mathcode`F="7046 \mathcode`G="7047 \mathcode`H="7048
\mathcode`I="7049 \mathcode`J="704A \mathcode`K="704B \mathcode`L="704C
\mathcode`M="704D \mathcode`N="704E \mathcode`O="704F \mathcode`P="7050
\mathcode`Q="7051 \mathcode`R="7052 \mathcode`S="7053 \mathcode`T="7054
\mathcode`U="7055 \mathcode`V="7056 \mathcode`W="7057 \mathcode`X="7058
\mathcode`Y="7059 \mathcode`Z="705A

% handling accented characters in plain TeX :

% \catcode`\À=\active \defÀ{\`A}    \catcode`\à=\active \defà{\`a} 
% \catcode`\Â=\active \defÂ{\^A}    \catcode`\â=\active \defâ{\^a} 
% \catcode`\Æ=\active \defÆ{\AE}    \catcode`\æ=\active \defæ{\ae}
% \catcode`\Ç=\active \defÇ{\c C}   \catcode`\ç=\active \defç{\c c}
% \catcode`\È=\active \defÈ{\`E}    \catcode`\è=\active \defè{\`e} 
% \catcode`\É=\active \defÉ{\'E}    \catcode`\é=\active \defé{\'e} 
% \catcode`\Ê=\active \defÊ{\^E}    \catcode`\ê=\active \defê{\^e} 
% \catcode`\Ë=\active \defË{\"E}    \catcode`\ë=\active \defë{\"e} 
% \catcode`\Î=\active \defÎ{\^I}    \catcode`\î=\active \defî{\^\i}
% \catcode`\Ï=\active \defÏ{\"I}    \catcode`\ï=\active \defï{\"\i}
% \catcode`\Ô=\active \defÔ{\^O}    \catcode`\ô=\active \defô{\^o} 
% \catcode`\Ù=\active \defÙ{\`U}    \catcode`\ù=\active \defù{\`u} 
% \catcode`\Û=\active \defÛ{\^U}    \catcode`\û=\active \defû{\^u} 
% \catcode`\Ü=\active \defÜ{\"U}    \catcode`\ü=\active \defü{\"u} 

\def\mod{\mathop{\rm mod.}\nolimits}
\def\pmod#1{\;(\mod#1)}

\openup1\jot

\font\tencyr=wncyr10
\font\sevencyr=wncyr7
\font\fivecyr=wncyr5
\newfam\cyrfam
\def\cyr{\fam\cyrfam\tencyr}
\textfont\cyrfam=\tencyr
\scriptfont\cyrfam=\sevencyr
\scriptscriptfont\cyrfam=\fivecyr

\newcount\refno

\long\def\ref#1:#2<#3>{                                        
\global\advance\refno by1\par\noindent                              
\llap{[{\bf\number\refno}]\ }{#1} {\it #2} #3\goodbreak }

% \long\def\ref#1:#2<#3>{                                        
% \global\advance\refno by1\par\noindent                              
% \llap{[{\bf\number\refno}]\ }{#1}\pointir{\it #2} #3\goodbreak }

\def\citer#1(#2){[{\bf\number#1}\if#2\empty\relax\else,\ #2\fi]}

\newbox\bibbox
\setbox\bibbox\vbox{
%\bigskip
%\centerline{---$*$---$*$---}
\bigbreak\bigbreak
\centerline{{\pc BIBLIOGRAPHY}}

\def\mr#1.{{\sevenrm MR#1.}}

\refno=18

\ref{\pc BARNET-\pc LAMB} (T), {\pc GERAGHTY} (D), {\pc  HARRIS}
 (M) \& {\pc TAYLOR} (R), :
A family of Calabi-Yau varieties and potential automorphy II, <Publications
of the Research Institute for Mathematical Sciences {\bf 47} (2011) 1, 29--98.>
\newcount\bght  \global\bght=\refno

\ref{\pc GEE} (T), :
%Elliptic curves,
<youtu.be/6eZQu120A80.>
\newcount\gee  \global\gee=\refno

\ref{\pc TAYLOR} (R), :
<youtu.be/EWHiMxfummQ>
\newcount\taylorclay \global\taylorclay=\refno

\ref{\pc BIRCH} (B) \& {\pc STEPHENS} (N), :
The parity of the rank of the Mordell-Weil group,
<Topology {\bf 5} (1966), 295--299.>
\newcount\birchstephens \global\birchstephens=\refno

\ref{\pc HEEGNER} (K), :
Diophantische analysis und modulfunktionen,
<Mathematische Zeitschrift {\bf 56} (1952), 227--253.>
\newcount\heegner  \global\heegner=\refno

\ref{\pc TIAN} (Y), :
Congruent numbers and Heegner points,
<Cambridge Journal of Mathematics {\bf 2} (2014) 1, 117--161.>
\newcount\tian  \global\tian=\refno

\ref {\pc COATES} (J), :
Congruent numbers,
<Acta Mathematica Vietnamica {\bf 39} (2014) 1, 3--10.>
\newcount\coates \global\coates=\refno

\ref{\pc ZHANG} (S), :
Congruent numbers and Heegner points,
<Asia Pacific Mathematics Newsletter {\bf 3} (2013) 2, 12--15.>
\newcount\zhang  \global\zhang=\refno

\ref{\pc ZHANG} (S), :
%The ongruent number problem and the Birch \& Swinnerton-Dyer conjecture,
% <youtu.be/\_ZYtuXEjdc8.>
<carmin.tv/s/6433.>
\newcount\luminy  \global\luminy=\refno

\ref{\pc SMITH} (A), :
%$2^\infty$-Selmer groups, $2^\infty$-class groups, and Goldfeld's conjecture,
<arXiv\string:1702.02325.>
\newcount\smith  \global\smith=\refno

\ref{\pc TIAN} (Y), {\pc YUAN} (X), \& {\pc ZHANG} (S), :
Genus periods, genus points and congruent number problem,
<Asian Journal of Mathematics {\bf 21} (2017) 4, 721--773.>
\newcount\yetianzhang \global\yetianzhang=\refno

\ref{\pc TIAN} (Y), :
Congruent number problem,
<Notices of the International Congress of
Chinese Mathematicians {\bf 5} (2017) 1, 51--61.>
\newcount\tianiccm \global\tianiccm=\refno

\ref{\pc KRIZ} (D), :
% Supersingular main conjectures, Goldfelds' conjecture and Silvesters'
% conjecture, 
% <Preprint (2020).>
<arXiv\string:2002.04767>
\newcount\kriz \global\kriz=\refno

\ref{\pc SMITH} (A), :
%$2^\infty$-Selmer groups, $2^\infty$-class groups, and Goldfeld's conjecture,
<youtu.be/mdZtLoTT1yk>
\newcount\ias  \global\ias=\refno

\ref{\pc SMITH} (A), :
%$2^k$-Selmer groups and Goldfeld's conjecture,
<carmin.tv/s/5462.>
\newcount\ihp  \global\ihp=\refno

% \ref{\pc COATES} (J) :
% <newton.ac.uk/seminar/24765/.> See also {\it The oldest problem},
% Notices of the International Congress of Chinese Mathematicians {\bf
% 5} (2017) 2, 8--13.
% \newcount\coatesiccm \global\coatesiccm=\refno

\ref{\pc COATES} (J), :
The oldest problem,
<Notices of the International Congress of Chinese Mathematicians {\bf
5} (2017) 2, 8--13. newton.ac.uk/seminar/24765/>  
\newcount\coatesiccm \global\coatesiccm=\refno

\ref{\pc BURUNGALE} (A) \& {\pc TIAN} (Y), :
The even parity Goldfeld conjecture: Congruent number elliptic curves,
<Journal of Number Theory {\bf 230} (2022), 161--195.>
\newcount\burungale  \global\burungale=\refno

% \ref{\pc LOZANO-\pc ROBLEDO} (A):
% <youtu.be/sHTACMV3IEc.>
% \newcount\alvaro \global\alvaro=\refno

\ref{\pc KATZ} (N) \& {\pc SARNAK} (P), :
Zeroes of zeta functions and symmetry,
<Bulletin of the American Mathematical Society {\bf 36} (1999) 1, 1--26.>
\newcount\katzsarnak  \global\katzsarnak=\refno

\ref{\pc BHARGAVA} (M) \& {\pc SHANKAR} (A), :
Binary quartic forms having bounded invariants, and the boundedness of
the average rank of elliptic curves,
<Annals of Mathematics {\bf 181} (2015) 1, 191--242.>
\newcount\bhargavai \global\bhargavai=\refno

\ref{\pc BHARGAVA} (M) \& {\pc SHANKAR} (A), :
Ternary cubic forms having bounded invariants, and the existence of a
positive proportion of elliptic curves having rank~$0$,
<Annals of Mathematics {\bf 181} (2015) 2, 587--621.>
\newcount\bhargavaii \global\bhargavaii=\refno

\ref{\pc BHARGAVA} (M), :
<youtu.be/t5S8PO2PQgU>
\newcount\bhargavavid \global\bhargavavid=\refno

\ref{\pc BHARGAVA} (M) \& {\pc SKINNER} (C), :
A positive proportion of elliptic curves over $\bf Q$ have rank~$1$,
<Journal of the Ramanujan Mathematical Society {\bf 29} (2014) 2,
221--242.>
\newcount\bhaskinner \global\bhaskinner=\refno

\ref{\pc BHARGAVA} (M), {\pc SKINNER} (C) \& {\pc ZHANG} (W), :
<arXiv\string:1407.1826.>
\newcount\bsz \global\bsz=\refno

\ref{\pc BHARGAVA} (M), :
<youtu.be/00WUkXaDm-s.>
\newcount\mbvid \global\mbvid=\refno

\ref{\pc POONEN} (B), :
Average rank of elliptic curves [after Manjul Bhargava \& Arul
Shankar],
<S\'eminaire Bourbaki, Ast\'erisque {\bf 352}
(2013), Expos\'e 1049.>
\newcount\poonenbou \global\poonenbou=\refno

\ref{\pc WILES} (A), :
<youtu.be/1WYlP-B9nPI.>
\newcount\wiles \global\wiles=\refno

\ref{\pc SKINNER} (C), :
<heilbronn.ac.uk/video-archive/>
\newcount\skinner \global\skinner=\refno

\ref{\pc GROSS} (B), :
<youtu.be/EXqRiaZYeNI.>
\newcount\gross \global\gross=\refno

\ref{\pc WILES} (A), :
<claymath.org/sites/default/files/birchswin.pdf.>
\newcount\wilesbsd \global\wilesbsd=\refno

\ref{\pc COATES} (J), :
Lectures on the Birch \& Swinnerton-Dyer conjecture,
<Notices of the International Congress of
Chinese Mathematicians {\bf 1} (2013) 2, 29--46.>
\newcount\coatesiccmii \global\coatesiccmii=\refno

 \ref{\pc ZHANG} (W), :
The Birch \& Swinnerton-Dyer conjecture and Heegner points: a survey,
<Current developments in mathematics, International Press, Somerville,
2014.>
\newcount\zhangwei \global\zhangwei=\refno

\ref{\pc BURUNGALE} (A), {\pc SKINNER} (C), \& {\pc TIAN} (Y), :
The Birch \& Swinnerton-Dyer conjecture: a brief survey,
<Proceeding of Symposposia in Pure Mathematics {\bf 104},
Amererican Mathematical Society, Providence, 2021.>
\newcount\buruskitian \global\buruskitian=\refno

\ref{\pc BURUNGALE} (A), :
<youtu.be/rVq-zFGx798.>
\newcount\burias \global\burias=\refno

\ref{\pc TIAN} (Y), :
<youtu.be/Ed6WNjXVNuc.>
\newcount\tianprinc \global\tianprinc=\refno

\ref{\pc BURUNGALE} (A), {\pc SKINNER} (C), \& {\pc TIAN} (Y), :
Elliptic curves and Beilinson-Kato elements~: rank one aspects,
<Preprint.>
\newcount\buruskitianpp \global\buruskitianpp=\refno

\ref{\pc KIM} (C-H), :
<arXiv\string:2109.12344.>
\newcount\chanhokim \global\chanhokim=\refno

\ref{\pc TATE} (J), :
<youtu.be/RtiVaALdqX0.>
\newcount\tate \global\tate=\refno

\ref{\pc COLMEZ} (P), :
La conjecture de Birch et Swinnerton-Dyer $p$-adique,
<S\'eminaire Bourbaki, Ast\'erisque {\bf 294} (2004), Expos\'e 919.>
\newcount \colmez \global\colmez=\refno

\ref{\pc MONSKY} (P), :
Generalizing the Birch-Stephens theorem, I: Modular curves,
<Mathematische Zeitschrift {\bf 221} (1996), 415–420.>
\newcount\monsky \global\monsky=\refno

\ref{\pc DOKCHITSER} (T), :
<arXiv\string:1009.5389.>
\newcount\timdok \global\timdok=\refno

\ref{\pc POONEN} (B), :
Heuristics for the arithmetic of elliptic curves,
<Proceedings of the International Congress of Mathematicians --- Rio
de Janeiro, World Scientific Publishers, Hackensack, 2018.>
\newcount\poonenrio \global\poonenrio=\refno

\ref{\pc POONEN} (B), :
<youtu.be/yxyz2K1UX9Y.>
\newcount\poonenriovid \global\poonenriovid=\refno

\ref{\pc POONEN} (B), :
<carmin.tv/s/5784>
\newcount\poonenriolum \global\poonenriolum=\refno

\ref{\pc ELKIES} (N), :
<youtu.be/xkznwW6jFh0.>
\newcount\elkies \global\elkies=\refno

}

\centerline{\bf Congruent numbers, elliptic curves,}
\centerline{\bf and the passage from the local to the
global} 
\centerline{\it An update}
\medskip
\centerline{Chandan Singh Dalawat}
\bigskip

\rightline{\it Number-theory is not standing still.}
\rightline{--- Andr{\'e} Weil, {\it Two lectures on Number Theory}, 1974}

\bigskip

Since the appearance of my popular article on {\it Congruent numbers,
elliptic curves, and the passage from the local to the global},
published in {\it Resonance\/} in December~2009, a number of new
results have been obtained on the topics discussed there.  We review a
few of these, using the same notation and terminology.  We also
continue with the numbering of the statements and bibliographic items.

Three major advances in the arithmetic of elliptic curves have taken
place in the last decade.  The Goldfeld Conjecture (recalled below)
has been proved for the family of elliptic curves arising from the
congruent number problem (and for many other families besides), a
converse to the theorems of Coates \& Wiles, Rubin, Gross \& Zagier,
Kolyvagin, and Kato has been proved, and a number of striking
statistical results has been obtained about elliptic curves whose
group of rational points has rank~0 or~1, about the average rank of
elliptic curves over $\bf Q$, and about the proportion which satisfy
the Birch \& Swinnerton-Dyer Conjecture.  There has also been much
speculation about the boundedness of the rank in various families of
elliptic curves, or about the number of elliptic curves of given rank.
These are the four topics we touch upon in this update.  We also
briefly mention that the Parity Conjecture (a very special case of the
Birch \& Swinnerton-Dyer Conjecture) follows from the Shafarevich-Tate
conjecture.

Before we begin, I should warn the reader that some of the material
covered in these notes has not yet appeared in print and therefore has not
yet been fully verified by the mathematical community.  My role here is
that of a reporter, not an endorser.  

First of all, note that the technical hypothesis in Theorem~35
requiring the elliptic curve to have at least one place of
multiplicative reduction has been removed~:

\th THEOREM 36 (T. Barnet-Lamb, D. Geraghty, M. Harris, \&
R. Taylor \citer\bght()) 
\enonce
The Sato-Tate conjecture (Conjecture~$34$) is true for every elliptic
curve over ${\bf Q}$ (without complex multiplications).
\endth

An elementary introduction to the problem of congruent numbers and
elliptic curves leading up to the Sato-Tate conjecture can be found in
the inaugural lecture by Toby Gee at the Imperial
College \citer\gee().  See also Richard Taylor's talk at the Clay
conference in 2007 \citer\taylorclay().

\bigbreak
{\bf The sign of the functional equation}

We have seen that a squarefree integer $\alpha>0$ is congruent if and
only if there exist $x,y\in{\bf Q}$ such that $y\neq0$ and $\alpha
y^2=x^3-x$ (Proposition~6).  In other words, the group $C_\alpha({\bf
Q})$ of rational points on the elliptic curve
$$
C_\alpha:\alpha y^2=x^3-x
$$
is infinite if and only if $\alpha$ is a congruent number.  The
conjecture of Birch \& Swinnerton-Dyer (Conjecture~24) predicts that
the group $C_\alpha({\bf Q})$ is infinite if and only if
$L(C_\alpha,1)=0$ (the $L$-function $L(C_\alpha,s)$ vanishes at $s=1$,
or equivalently the order of vanishing ${\rm ord}_{s=1}L(C_\alpha,s)$
is $>0$).  By computing the ``global root number'' or the ``sign of
the functional equation'' of $C_\alpha$ \citer\birchstephens(), one
can check that ${\rm ord}_{s=1}L(C_\alpha,s)$ is {\it even\/} if
$\alpha\equiv1,2,3\pmod8$ and {\it odd\/} if
$\alpha\equiv-3,-2,-1\pmod8$.

If the order of vanishing is odd, then certainly $L(C_\alpha,1)=0$.
We are thus led by Birch \& Swinnerton-Dyer to the following~:

\th CONJECTURE 37
\enonce
If\/ $\alpha>0$ is squarefree and\/ $\equiv-3,-2,-1\pmod8$, then\/
it is a congruent number.
\endth

The converse is not true~: 34 is the area of the right-triangle with
sides $225/30$, $272/20$, and $353/30$.  Nearly a decade before
Birch \& Swinnerton-Dyer came up with their conjecture, Heegner had
already taken the first important step towards proving Conjecture~37
in a paper \citer\heegner() which went unnoticed at the time~:

\th THEOREM 38 (K. Heegner, 1952)
\enonce
Let $\alpha>0$ be a  squarefree integer such that\/ 
$\alpha\equiv-3,-2,-1\pmod8$.  If\/ $\alpha$ has precisely one odd prime 
divisor, then\/ $\alpha$ is a congruent number.
\endth

This was generalised by Monsky (1990) to allow $\alpha$ to have two
odd prime divisors.  Ye Tian \citer\tian() took a major step forward
by proving the following important result~:

\th THEOREM 39 (Y. Tian, 2014)
\enonce
For every $j\equiv-3,-2,-1\pmod8$ and every $n>0$, there are infinitely
many squarefree integers\/ $\alpha>0$ having precisely\/ $n$ odd prime
divisors and such that\/ $\alpha\equiv j\pmod8$ which are congruent
numbers.
\endth

For an introduction to Tian's work, see the article \citer\coates() by
John Coates.

How about squarefree integers $\alpha>0$ such that
$\alpha\equiv1,2,3\pmod8$, when the order of vanishing ${\rm
ord}_{s=1}L(C_\alpha,s)$ is {\it even\/}~?  If moreover $\alpha$ has a
unique odd prime divisor~$p$, there is the following result which goes
back to the 19$^{\rm th}$ century in part~:

\th THEOREM 40 (A. Genochhi, 1874~; M. Razar, 1974)
\enonce
Let\/ $p$ be a prime number. If\/ $p\equiv3\pmod8$, then\/ $p$ is not
congruent.  If\/ $p\equiv5\pmod8$, then\/ $2p$ is not congruent.   
\endth

What happens if $p\equiv1\pmod8$, or if $\alpha$ has several odd prime
divisors~?  The following analogue of Theorem~39 says something about
these cases~:

\th THEOREM 41 (K. Feng, 1996~; D. Li \& Y. Tian, 2000~; C. Zhao,
2001)
\enonce
For every $j\equiv1,2,3\pmod8$ and every $n>0$, there are infinitely
many squarefree integers\/ $\alpha>0$ having precisely\/ $n$ odd prime
divisors and such that\/ $\alpha\equiv j\pmod8$ which are not congruent
numbers.
\endth

See for example the expository article \citer\zhang() by Shou-Wu Zhang
and his lecture \citer\luminy() at Luminy.

\bigbreak
{\bf The minimalist philosophy}

Before we state Goldfeld's conjecture in the special case of the
family of elliptic curves $C_\alpha$ (where $\alpha>0$ is a squarefree
integer), we need to recall that if $S$ is a set of positive integers
and $T$ is a subset of $S$, then the {\it density\/}, or the {\it
natural density\/}, of $T$ in $S$ is the limit (if it exists)
$$
\delta_{T,S}=\lim_{x\to+\infty}{t_x\over s_x},
$$
where $t_x$ (resp.  ~$s_x$) is the number of elements in $T$
(resp.~$S$) which are $<x$.  We also say that $100\delta_{T,S}$\% of
elements of $S$ are in $T$. For example, when $S={\bf N}$ is the set
of all positive integers and $T=2{\bf N}$ is the set of even numbers,
then $\delta_{S,T}=1/2$, or $50$\% of positive integers are even.  We
also say that the {\it proportion\/} of elements of $S$ which are in
$T$ is $\delta_{T,S}$.

Note that $\delta_{T,S}=0$ does not imply that $T$ is empty, nor does
$\delta_{T,S}=1$ imply that $T=S$.  For example, 0\% of positive
integers are prime numbers, yet there are infinitely many of them.

\th CONJECTURE 42 (D. Goldfeld, 1979)
\enonce
Among the squarefree integers $\alpha>0$ such that
$\alpha\equiv1,2,3\pmod8$, the percentage of those for which\/ ${\rm
ord}_{s=1}L(C_\alpha,s)=0$
% ${\rm rk}_{\bf Z}\,C_\alpha({\bf Q})=0$
is\/~$100$\%. 
\endth

This is sometimes called Goldfeld's even-parity conjecture.  Combined
with the theorem of Coates \& Wiles, this means that the percentage of
such $\alpha$ for which $C_\alpha({\bf Q})$ is infinite should be 0\%.
As $\alpha$ is a congruent number if and only if $C_\alpha({\bf Q})$
is infinite, Conjecture~42 implies the following~

\th CONJECTURE 43 (D. Goldfeld, 1979)
\enonce
The percentage of squarefree $\alpha>0$ such that
$\alpha\equiv1,2,3\pmod8$ which are congruent is\/~$0$\%.
\endth

Similary, he conjectured the following odd-parity version~:

\th CONJECTURE 44 (D. Goldfeld, 1979)
\enonce
Among the squarefree integers $\alpha>0$ such that
$\alpha\equiv-3,-2,-1\pmod8$, the percentage of those for which\/ ${\rm
ord}_{s=1}L(C_\alpha,s)=1$
% ${\rm rk}_{\bf Z}\,C_\alpha({\bf Q})=0$
is\/~$100$\%.
\endth

By Kolyvagin's theorem (Theorem 28), this would imply that the
percentage of squarefree $\alpha>0$ such that
$\alpha\equiv-3,-2,-1\pmod8$ which are congruent numbers is 100\%.

Alexander Smith \citer\smith() and Daniel Kriz \citer\kriz() have
succeeded in proving these three conjectures but the proof has not yet
been published.  In fact, their result applies to the familiy of
``quadratic twists'' of any elliptic curve satisfying a certain
technical hypothesis (which $C_1$ does).

\th THEOREM 45 (A. Smith, 2017)
\enonce
The percentage of squarefree $\alpha>0$ such that
$\alpha\equiv1,2,3\pmod8$ which are congruent is\/~$0$\%.
\endth

% Using the work of Ye Tian, Xinyi Yuan and Shou-Wu
% Zhang \citer\yetianzhang() (see also \citer\tianiccm() for an
% expository account) on the one hand, and the theorem of Gross \&
% Zagier on the other, he also obtained the following result in the
% odd-parity case~:

% \th THEOREM 44 (A. Smith, 2017)
% \enonce
% At least $50$\% of squarefree $\alpha>0$ such that
% $\alpha\equiv-3,-2,-1\pmod8$ are congruent numbers.  If the Birch \&
% Swinnerton-Dyer conjecture is true, then the percentage is\/ $100$\%.
% \endth

One can listen to Smith's talks at the Instiute for Advanced
Study \citer\ias() and at the Institut Henri
Poincar{\'e} \citer\ihp().  For an expository account,
see \citer\coatesiccm().  

\th THEOREM 46 (D. Kriz, 2020)
\enonce
The percentage of squarefree $\alpha>0$ such that
$\alpha\equiv-3,-2,-1\pmod8$ which are congruent is\/~$100$\%.
\endth

Note that the Birch \& Swinnerton-Dyer Conjecture predicts that {\it
all\/} squarefree $\alpha>0$ such that $\alpha\equiv-3,-2,-1\pmod8$
are congruent, because ${\rm ord}_{s=0}L(C_\alpha,s)$ is odd, as
explained above (Conjecture~37).

Recall (Notation 10) that we had defined $c_j(n)$, for\/ $j=1,2$ and
$n>0$, as the coefficient of\/ $T^n$ in the formal series \/
$g(T)\theta_j(T)$, where
$$
g(T)=T\prod_{n>0}(1-T^{8n})(1-T^{16n})\quad\hbox{and}\quad
\theta_j(T)=1+2\sum_{n>0}T^{2jn^2}. 
$$

Recall also that if $\alpha=jn$ with $j=1,2$ and $n>0$ squarefree and
{\it odd\/}, then $\alpha$ is not a congruent number if $c_j(n)\neq0$,
by the theorems of Coates \& Wiles (Theorem~14) and Tunnell
(Thorem~25).  In view of this observation, Theorem~45 follows from the
following recent result~:

\th THEOREM 47 (A. Burungale \& Y. Tian \citer\burungale())
\enonce
The percentage of squarefree integers $\alpha=jn$ (with $j=1,2$ and
$n>0$ odd) such that $\alpha\equiv1,2,3\pmod8$ for which $c_j(n)=0$
is\/~$0$\%.
\endth

% For an overview of some of these questions, see the interesting online
% talk by Alvaro Lozano-Robledo \citer\alvaro().

\bigbreak
{\bf Statistics of elliptic curves}

We have just come across conjectures about the percentage of elliptic
curves in the family $C_\alpha$ ($\alpha>0$ running through squarefree
integers) for which the group $C_\alpha({\bf Q})$ has rank~$0$,
rank~$1$, etc.  The curves $C_\alpha$ are ``quadratic twists'' of the
elliptic curve $C_1$.  There are similar conjectures about the the
family of all elliptic curves over $\bf Q$, giving the percentage of
curves for which the group of rational points has rank~$0$, rank~$1$,
etc.  These percentages are defined by first computing the percentage,
among the finitely many elliptic curves of bounded ``height'', of
those which have the given rank, and then taking the limit as the
bound goes to $+\infty$.  More precisely, every elliptic curve $E$
over $\bf Q$ admits a unique equation
$$
E:\quad y^2=x^3+A_4x+A_6,\qquad(A_4,A_6\in{\bf Z}~;\ \ 4A_4^3+27A_6^2\neq0)
$$
such that {\it there is no prime number\/~$p$ for which\/ $p^{12}$
divides ${\rm gcd}(A_4^3,A_6^2)$}~; the {\it height\/} of $E$ is defined
to be $h(E)={\rm Sup}(|4A_4^3|, |27A_6^2|)$. Clearly, for every $c>0$,
there are only finitely many elliptic curves $E$ over $\bf Q$ such
that $h(E)<c$.

One can similarly define the {\it average rank of elliptic curves
over\/} $\bf Q$ as the limit, when $c\to+\infty$, of the average of
${\rm rk}_{\bf Z}E({\bf Q})$ for the finitely many elliptic curves $E$
over $\bf Q$ such that $h(E)<c$.

\th CONJECTURE 48 (D. Goldfeld, 1979~; N. Katz \&
P. Sarnak \citer\katzsarnak())
\enonce
For $50$\% of elliptic curves\/ $E$ over\/ $\bf Q$, the group\/
$E({\bf Q})$ has rank\/~$0$~; for\/ $50$\% of them, it has
rank\/~$1$~; for\/ $0$\% of them, it has rank\/~$>1$.  Consequently,
the average rank of elliptic curves over\/ $\bf Q$ is ${1\over 2}$.
\endth

When this conjecture was made, it was not known that the precentage of
elliptic curves of rank $0$ or $1$ is $>0$\%, let alone $=100$\%, or
that the average rank of elliptic curves is $<+\infty$, let alone
$={1\over2}$.  Since then, remarkable progress has been made~:

\th THEOREM 49 (M. Bhargava \&
A. Shankar \citer\bhargavai() )
\enonce
When elliptic curves\/ $E$ over\/ $\bf Q$ are ordered by height, the
average rank of the\/ $\bf Z$-module\/ $E({\bf Q})$ is\/ $<1$.
\endth

\th THEOREM 50 (M. Bhargava \&
A. Shankar \citer\bhargavaii() )
\enonce
When elliptic curves\/ $E$ over\/ $\bf Q$ are ordered by height, the
proportion for which\/ ${\rm rk}_{\bf Z}E({\bf Q})=0$ is\/ $>0$, and
the proportion which satisfy the Birch \& Swinnerton-Dyer conjecture
is\/ $>0$.
\endth

\th THEOREM 51 (M. Bhargava \& C. Skinner \citer\bhaskinner() )
\enonce
When elliptic curves\/ $E$ over\/ $\bf Q$ are ordered by height, the
proportion for which\/ ${\rm rk}_{\bf Z}E({\bf Q})=1$ is\/ $>0$.
\endth

For a succint account of these results when they were discovered, see
Manjul Bhargava's talk at the Clay conference at Harvard in
2011 \citer\bhargavavid()~; see also the
Bourbaki talk by Bjorn Poonen \citer\poonenbou().  Since then the percentages have been
steadily improved.  For example, within a few years, it was proved
that a majority of elliptic curves satisfy Conjecture~24~:

\th THEOREM 52 (M. Bhargava, C. Skinner, \& W. Zhang \citer\bsz() )
\enonce
When ordered by height, more than half the elliptic curves over\/ $\bf
Q$ satisfy the Birch \& Swinnerton-Dyer conjecture.
\endth

For recent oral overviews of the Birch \& Swinnerton-Dyer conjecture,
see the talks by Manjul Bhargava \citer\mbvid() and Andrew
Wiles \citer\wiles() at Clay conferences, by Chris Skinner at the
Heilbronn Institute \citer\skinner(), and by Benedict
Gross \citer\gross() at Caltech.  Among the written sources, one may
consult the description by Wiles \citer\wilesbsd() for the Clay
Institute, the lectures by John Coates \citer\coatesiccmii(), the
article by Wei Zhang \citer\zhangwei(), and the recent article by
Ashay Burungale, Chris Skinner \& Ye Tian \citer\buruskitian().

Let us finally mention that there has been very important work towards
proving a converse to the theorems of Coates-Wiles (cf.~Thorem~14),
Rubin (cf.~Theorem~15), Gross \& Zagier, Kolyvagin (Theorem~28), and
Kato by Chris Skinner, Wei Zhang, Ashay Burungale, Ye Tian, Daniel
Kriz, and Chan-Ho Kim.  The proof of the following converse theorem is
now complete, although all the details have not yet been published~:

\th THEOREM 53 (A. Burungale, C. Skinner, \& Y. Tian \citer\buruskitianpp()) 
\enonce
Let\/ $E$ be an elliptic curve over\/ $\bf Q$. If the Shafarevich-Tate
conjecture (Conjecture~$26$) holds for $E$ and\/ ${\rm rk}_{\bf
Z}E({\bf Q})=0,1$, then\/ ${\rm ord}_{s=1}L(E,s)= {\rm rk}_{\bf
Z}E({\bf Q})$.
\endth

See also the preprint by Chan-Ho Kim \citer\chanhokim().
The interested reader is referred to the Princeton talks of
Burungale \citer\burias() and Tian \citer\tianprinc().  Note that
according to Conjecture~48, 100\% of elliptic curves over $\bf Q$ have
${\rm rk}_{\bf Z}E({\bf Q})=0,1$.

\bigbreak
{\bf The parity of the rank}

For an elliptic curve $E$ over ${\bf Q}$, the Parity Conjecture
asserts that ${\rm rk}_{\bf Z}E({\bf Q})\equiv{\rm
ord}_{s=1}L(E,s)\pmod2$.  Paul Monsky \citer\monsky() has proved that
this easy consequence of Conjecture~24 (Birch \& Swinnerton-Dyer)
follows from Conjecture~26 predicting the finiteness of the
Shafarevich-Tate group ${\cyr Sh}(E)$.  More recently, Tim \& Vladimir
Dokchitser have proved that it is enough to require the finiteness of
the $p$-primary component of ${\cyr Sh}(E)$ for a single arbitrary
prime $p$.  See the nice survey \citer\timdok().

\bigbreak
{\bf The boundedness of the rank}

For a long time it was believed that there are only finitely many
possibilities for the rank of the $\bf Z$-module $E({\bf Q})$ as $E$
runs through all elliptic curves defined over ${\bf Q}$.  Then it was
thought for a while that the ranks were unbounded, by analogy with
what happens over function fields over finite fields (Shafarevich \&
Tate, 1967~; Ulmer, 2002)~; see for example Tate's lecture on the
arithmetic of elliptic curves at Oslo \citer\tate() or the Bourbaki
talk by Colmez on the $p$-adic Birch \& Swinnerton-Dyer
conjecture \citer\colmez().  Noam Elkies has constructed infinitely
many $E$ for which the rank of $E({\bf Q})$ is $>18$, and one $E$ for
which ${\rm rk}_{\bf Z}E({\bf Q})>27$.

By constructing a linear algebraic model for the rank and other
algebraic invariants of elliptic curves, the consensus seems to have
returned to boundedness of late~:

\th CONJECTURE 54
\enonce
There are only finitely many elliptic curves defined over ${\bf Q}$
for which\/ ${\rm rk}_{\bf Z}E({\bf Q})>21$.
\endth

There is in fact a whole series of conjectures inspired by a
linear-algebraic model of elliptic curve invariants.  For example,
although there are infinitely many elliptic curves over $\bf Q$ such
that $E({\bf Q})$ has a point of order~$12$, only finitely many of
these should have rank $>2$, conjecturally.  See for example Bjorn
Poonen's lecture at Rio de Janeiro (text \citer\poonenrio(),
video \citer\poonenriovid()) and his presentation at
Luminy \citer\poonenriolum().  For more speculation on the rank, see
the online talk by Noam Elkies \citer\elkies().

\medbreak
{\bf Acknowledgements} I warmly thank Daniel Kriz and Ashay Burungale
for clarifying certain points, Xavier Xarles for correcting a link,
and Peter Scholze for his comments.

\bigbreak
\unvbox\bibbox
\vskip5cm plus1cm minus1cm
{\obeylines\parskip=0pt\parindent=0pt Chandan Singh Dalawat
Harish-Chandra Research Institute
Chhatnag Road, Jhunsi
{\pc ALLAHABAD} 211\thinspace019, India
\vskip5pt
\tt dalawat@gmail.com}

\bye